\documentclass[12pt,a4paper]{article}
\usepackage{amsmath,amsxtra,amsfonts,amssymb,latexsym,amsthm,amscd,verbatim,amstext}
\usepackage[mathscr]{eucal}
\newcommand{\me}{\mathrm{e}}

\newcommand{\dif}{\mathrm{d}}

\input cyracc.def

\newfont{\tricyr}{wncyr10 at 12pt}
\newfont{\tricyi}{wncyi10 at 12pt}
\newfont{\tricyb}{wncyb10 at 12pt}
\newfont{\Tricyr}{wncyr10 at 13.6pt}
\newfont{\Tricyi}{wncyi10 at 13.6pt}
\newfont{\Tricyb}{wncyb10 at 13.6pt}
\newfont{\tricmr}{cmr10 at 13.6pt}
\newfont{\tricmi}{cmti10 at 13.6pt}
\newfont{\tricmb}{cmb10 at 13.6pt}

\theoremstyle{plain}
\newtheorem{Th}{Theorem}
\newtheorem{Md}{Proposition}
\newtheorem{Lem}{Lemma}

\theoremstyle{definition}
\newtheorem{Def}{Definition}
\newtheorem{Not}{Remark}

\begin{document}
\centerline {{\bf Well-posedness for the Navier-Stokes equations with data}} 
\centerline {{\bf in homogeneous Sobolev-Lorentz spaces}} 

\vskip 0.7cm
\begin{center}
D. Q. Khai, N. M. Tri
\end{center}
\begin{center}
Institute of Mathematics, VAST\\
18 Hoang Quoc Viet, 10307 Cau Giay, Hanoi, Vietnam
\end{center}
\vskip 0.9cm 
{\bf Abstract}: In this paper, we study local well-posedness for the Navier-Stokes equations (NSE) 
with the arbitrary initial value in homogeneous Sobolev-Lorentz spaces $\dot{H}^s_{L^{q, r}}(\mathbb{R}^d):= (-\Delta)^{-s/2}L^{q,r}$
 for $d \geq 2, q > 1, s \geq 0$, $1 \leq r \leq \infty$, and $ \frac{d}{q}-1 \leq s < \frac{d}{q}$, this result improves the known 
results for $q > d,r=q, s = 0$ (see \cite{M. Cannone 1995, M. Cannone. Y. Meyer 1995})
 and for $q =r= 2, \frac{d}{2} - 1 < s  < \frac{d}{2}$ (see \cite{M. Cannone 1995, J. M. Chemin 1992}). \\
In the case of critical indexes ($s=\frac{d}{q}-1$), we prove global well-posedness for NSE provided the norm of the initial 
value is small enough. The result that is a generalization of the result in \cite{M. Cannone 1997} for $q = r=d, s=0$.


\footnotetext[3]{{\it e-mail address}: khaitoantin@gmail.com triminh@math.ac.vn}
\vskip 1cm 
\centerline{\S1. Introduction} 
\vskip 0,5cm 
We consider the Navier-Stokes equations in $\mathbb R^d$:  
\begin{align} 
\left\{\begin{array}{ll} \partial _tu  = 
\Delta u - \nabla .(u \otimes u) - \nabla p , & \\ 
\nabla .u = 0, & \\
u(0, x) = u_0, 
\end{array}\right . \notag
\end{align}
which is a condensed writing for
\begin{align} 
\left\{\begin{array}{ll} 1 \leq k \leq d, \ \  \partial _tu_k  
= \Delta u_k - \sum_{l =1}^{d}\partial_l(u_lu_k) - \partial_kp , & \\ 
\sum_{l =1}^{d}\partial_lu_l = 0, & \\
1 \leq k \leq d, \ \ u_k(0, x) = u_{0k} .
\end{array}\right . \notag
\end{align}
The unknown quantities are the velocity 
$u(t, x)=(u_1(t, x),\dots,u_d(t, x))$ of the fluid 
element at time $t$ and position $x$ and the pressure $p(t, x)$.\\
In the 1960s, mild solutions were first  constructed by 
Kato and Fujita (\cite{T. Kato 1962}, \cite{H. Fujita 1964}) 
that are continuous in time and take values in the Sobolev spaces 
$H^s(\mathbb{R}^d), (s \geq \frac{d}{2} - 1)$, say 
$u \in C([0, T]; H^s(\mathbb{R}^d))$. In 1992, a modern treatment 
for mild solutions in $H^s(\mathbb{R}^d), (s \geq \frac{d}{2} - 1)$ 
was given by Chemin \cite{J. M. Chemin 1992}. 
In 1995, using the simplified version of the bilinear operator, 
Cannone proved the existence of mild solutions in 
$\dot{H}^s(\mathbb{R}^d), (s \geq \frac{d}{2} - 1)$, 
see \cite{M. Cannone 1995}. 
Results on the existence of mild solutions with value in 
$L^q(\mathbb{R}^d), (q > d)$ were established in  the papers 
of  Fabes, Jones and Rivi\`{e}re \cite{E. Fabes 1972a} 
and of Giga \cite{Y. Giga 1986a}. Concerning the initial 
data in the space $L^\infty$, the existence of a mild solution was 
obtained by Cannone and Meyer in (\cite{M. Cannone 1995}, 
\cite{M. Cannone. Y. Meyer 1995}). 
In 1994, Kato and 
Ponce \cite{T. Kato 1994} showed that the NSE are well-posed 
when the initial data belong to the homogeneous Sobolev 
spaces $\dot{H}^{\frac{d}{q} - 1}_q(\mathbb{R}^d), (d \leq q < \infty)$. 
Recently, the authors of this article have 
considered NSE in the mixed-norm Sobolev-Lorentz spaces, 
see \cite{N. M. Tri: Tri2014a}.\\
In this paper,  for $d \geq 2, q > 1, s \geq 0$, $1 \leq r \leq \infty$, and $ \frac{d}{q}-1 \leq s < \frac{d}{q}$, we \linebreak investigate mild solutions to NSE in the spaces $L^\infty\big([0, T]; \dot{H}^s_{L^{q, r}}(\mathbb{R}^d)\big)$  when the initial data  belong  
to the Sobolev-Lorentz spaces $\dot{H}^s_{L^{q,r}}(\mathbb{R}^d)$, which are more general than the spaces 
$\dot{H}^s_q(\mathbb{R}^d)$, ($\dot{H}^s_q(\mathbb{R}^d) = \dot{H}^s_{L^{q,q}}(\mathbb{R}^d) $). 
We obtain the existence of mild solutions with arbitrary initial 
value when $T$ is small enough, and existence of mild solutions 
for any $T > 0$ when the norm of the initial value in 
the Besov spaces $\dot{B}^{s- d(\frac{1}{q} - \frac{1}{\tilde q}), 
\infty}_{\tilde q}(\mathbb R ^d)$, $\big(\frac{1}{2}
({\frac{1}{q}+\frac{s}{d}}) < \frac{1}{\tilde q} < 
{\rm min} \big\{\frac{1}{2} +\frac{s}{2d}, \frac{1}{q}\big\}\big)$ is small enough.\\
 In the particular case $(q > d,r=q, s = 0)$, we get the result which is 
more general than that of Cannone and Meyer (\cite{M. Cannone 1995}, 
\cite{M. Cannone. Y. Meyer 1995}). Here we obtained a statement  
that is stronger than that  of Cannone and Meyer but under a much weaker 
condition on the initial data.\\
In the particular case $(q =r= 2, \frac{d}{2} - 1 < s  < \frac{d}{2})$, 
we get the result which is more general than those of Chemin  
in \cite{J. M. Chemin 1992} and Cannone in \cite{M. Cannone 1995}. 
Here we obtained a statement that is stronger than those 
of Chemin  in \cite{J. M. Chemin 1992} and Cannone 
in \cite{M. Cannone 1995} but under a much weaker 
condition on the initial data.\\
In the case of critical indexes $(1 < q \leq d, r \geq 1, s = \frac{d}{q} - 1)$, 
we get a \linebreak result that is a generalization of a result of Cannone 
\cite{M. Cannone 1997}. In particular, when $q = r=d, s=0$, 
we get back the Cannone theorem (Theorem 1.1 in \cite{M. Cannone 1997}).\\
The paper is organized as follows. In Section $2$ we prove some inequalities 
for pointwise products in the Sobolev spaces and some auxiliary lemmas. 
In Section $3$ we present the main results of the paper.   
In the sequence, for a space of functions defined on $\mathbb R^d$, 
say $E(\mathbb R^d)$, we will abbreviate it as $E$. 
\vskip 0.5cm 
\centerline{\S2. Some auxiliary results}
\vskip 0,5cm
In this section, we recall the following results and notations.
\begin{Def}\label{def2.3} (Lorentz spaces). (See \cite{J. Bergh: 1976}.)\\
For $1 \leq p, r \leq \infty$, the Lorentz space $L^{p, r}(\mathbb{R}^d)$
is defined as follows: \linebreak A measurable function 
$f \in L^{p,r}(\mathbb{R}^d)$ if and only if\\
$\big\|f\big\|_{L^{p,r}}(\mathbb{R}^d) := 
\big(\int_0 ^\infty (t^{\frac{1}{p}}f^*(t))^r
\frac{\dif t}{t}\big)^{\frac{1}{r}} < \infty$ when $1 \leq r < \infty$,\\
$\big\|f\big\|_{L^{p,\infty}}(\mathbb{R}^d) :=  \underset{t > 0}
{\rm sup}\ t^{\frac{1}{p}}f^*(t) < \infty$ when $r = \infty $,\\
 where $f^*(t) = \inf \big\{\tau : \mathcal{M}^d (\{x: |f(x)| > \tau\}) \leq t \big\}$, 
with $\mathcal{M}^d$ being the  Lebesgue measure in $\mathbb R^d$.
\end{Def} 
Before proceeding to the definition of Sobolev-Lorentz spaces, 
let us introduce several necessary notations. For real number s, 
the operator $\dot{\Lambda}^s$ is defined through Fourier translation by 
$$
\big(\dot{\Lambda}^sf\big)^{\land}(\xi) = |\xi|^s\hat{f}(\xi).
$$
For $0 < s < d$, the operator $\dot{\Lambda}^s$ can be 
viewed as the inverse of the Riesz potential $I_s$ up to a positive constant
$$
I_s(f)(x) = \int_{\mathbb{R}^d}\frac{f(y)}{|x-y|^{d-s}}
\ \dif y\ \ \text{for}\  x \in \mathbb{R}^d.
$$
For $q>1, r \geq 1$, and $0 \leq s < \frac{d}{q}$, the operator $I_s$ 
is continuous from $L^{q, r}$ to $L^{\tilde q, r}$, where $\frac{1}{\tilde q} = \frac{1}{q}-\frac{s}{d}$, see (\cite{P. G. Lemarie-Rieusset 2002}, 
Theorem 2.4 $iii)$, p. 20).
\begin{Def}\label{def2.5}(Sobolev-Lorentz spaces). (See \cite{Hajaiej  Hichem}.)\\
For $q>1, r \geq 1$, and $0 \leq s < \frac{d}{q}$, the Sobolev-Lorentz space
$\dot{H}^s_{L^{q,r}}(\mathbb{R}^d)$ is defined 
as the space  $I_s(L^{q,r}(\mathbb{R}^d))$, equipped with the norm 
$$
\big\|f\big\|_{\dot{H}^s_{L^{q, r}}}: = \big\|\dot{\Lambda}^sf\big\|_{L^{q, r}}.
$$ 
\end{Def}
\begin{Lem}\label{lem2.1.1} Let $q > 1, 1 \leq r \leq \tilde{r}\leq \infty,
\  and\ 0 \leq s < \frac{d}{q}$. Then we have the following 
imbedding maps\\
{\rm (a)} 
\begin{gather*}
\dot{H}^s_{L^{q,1}} \hookrightarrow 
\dot{H}^s_{L^{q,r}} \hookrightarrow 
\dot{H}^s_{L^{q,\tilde{r}}} 
\hookrightarrow \dot{H}^s_{L^{q,\infty}}.
\end{gather*}
{\rm (b)} $\dot{H}^s_q = 
\dot{H}^s_{L^{q,q}}$ (equality of the norm).
\end{Lem}
\textbf{Proof}. It is easily deduced from the 
properties of the standard Lorentz spaces. \qed\\
In the following lemmas, we estimate the pointwise product of two 
functions in $\dot{H}^s_q(\mathbb{R}^d), (d \geq 2)$ 
which is a generalization of the Holder inequality. 
In the case when $s = 0$ we get back the usual  Holder inequality. 
Pointwise multiplication results for Sobolev spaces are also 
obtained in literature, see for example \cite{J.Y. Chemin}, 
\cite{P. G. Lemarie-Rieusset 2002}, \cite{T. Kato} 
and the references therein. 

\begin{Lem}\label{lem2.1}
Assume that 
$$
1 < p, q < d,\  and\  \frac{1}{p} + \frac{1}{q} < 1 +  \frac{1}{d}.
$$
Then the following inequality holds
$$
\big\|uv\big\|_{\dot{H}^1_r} \lesssim 
\big\|u\big\|_{\dot{H}^1_p}\big\|v\big\|_{\dot{H}^1_q}, 
\ \forall  u \in \dot{H}^1_p, v \in \dot{H}^1_q,
$$
where $\frac{1}{r} = \frac{1}{p} + \frac{1}{q} - \frac{1}{d}$.
\end{Lem}
\textbf{Proof}. 
By applying the Leibniz formula for the derivatives 
of a product of two functions, we have 
\begin{gather*}
 \big\|uv\big\|_{\dot{H}^1_r} \simeq \sum_{|\alpha| = 1}
\big\|\partial^{\alpha}(uv)\big\|_{L^r} \leq \sum_{|\alpha| = 1}
\big\|(\partial^{\alpha}u)v\big\|_{L^r} + \sum_{|\alpha| = 1}
\big\|u(\partial^{\alpha}v)\big\|_{L^r}.
\end{gather*}
By applying the H\"{o}lder and Sobolev inequalities we obtain 
\begin{gather*}
\sum_{|\alpha| = 1}\big\|(\partial^{\alpha}u)v\big\|_{L^r} 
\leq \sum_{|\alpha| = 1}\big\|\partial^{\alpha}u\big\|_{L^p}
\big\|v\big\|_{L^{q_1}} \lesssim \big\|u\big\|_{\dot{H}^1_p}
\big\|v\big\|_{\dot{H}^1_q},
\end{gather*}
where
$$
\frac{1}{q_1} = \frac{1}{q} - \frac{1}{d}.
$$
Similar to the above reasoning, we have
\begin{gather*}
\sum_{|\alpha| = 1}\big\|u(\partial^{\alpha}v)
\big\|_{L^r}  \lesssim \big\|u\big\|_{\dot{H}^1_p}
\big\|v\big\|_{\dot{H}^1_q}.
\end{gather*}
This gives the desired result
$$
\big\|uv\big\|_{\dot{H}^1_r} \lesssim 
\big\|u\big\|_{\dot{H}^1_p}\big\|v\big\|_{\dot{H}^1_q}.
$$
\qed
\begin{Lem}\label{lem2.1'}
Assume that 
\begin{equation}\label{them1}
0 \leq s \leq 1 , \frac{1}{p} > \frac{s}{d}, \frac{1}{q}
 > \frac{s}{d},\  and\  \frac{1}{p} + \frac{1}{q} < 1 +  \frac{s}{d}.
\end{equation}
Then the following inequality holds
$$
\big\|uv\big\|_{\dot{H}^s_r} \lesssim 
\big\|u\big\|_{\dot{H}^s_p}\big\|v\big\|_{\dot{H}^s_q}, 
\ \forall  u \in \dot{H}^s_p, v \in \dot{H}^s_q,
$$
where $\frac{1}{r} = \frac{1}{p} + \frac{1}{q} - \frac{s}{d}$.
\end{Lem}
\textbf{Proof}. It is not difficult to show that if 
$p, q, {\rm and}\ s$ satisfy \eqref{them1} then 
there exists numbers $p_1,p_2,q_1,q_2\in(1,+\infty)$ 
(may be many of them) such that 
\begin{gather*}
\frac{1}{p} = \frac{1 - s}{p_1} + \frac{s}{p_2},  \frac{1}{q} 
= \frac{1 - s}{q_1} + \frac{s}{q_2}, \frac{1}{p_1} + \frac{1}{q_1} < 1, \\ 
 p_2 < d, q_2 < d,\ {\rm and}\ \frac{1}{p_2} + \frac{1}{q_2} < 1 +  \frac{1}{d}.
\end{gather*}
Setting
$$
\frac{1}{r_1} = \frac{1}{p_1} + \frac{1}{q_1}, 
\frac{1}{r_2} = \frac{1}{p_2} + \frac{1}{q_2} - \frac{1}{d}, 
$$
we have
$$
\frac{1}{r} = \frac{1 - s}{r_1} + \frac{s}{r_2}.
$$
Therefore, applying Theorem 6.4.5 (page 152) of \cite{J. Bergh: 1976} 
(see also \cite{N. Kalton: 2007a} for $\dot H_{p}^{s}$), we get
\begin{gather*}
\dot{H}^s_p = [L^{p_1}, \dot{H}^1_{p_2}]_{s},  \dot{H}^s_q 
= [L^{q_1}, \dot{H}^1_{q_2}]_{s}, \dot{H}^s_r = [L^{r_1}, 
\dot{H}^1_{r_2}]_{s}.  
\end{gather*}
Applying the Holder inequality and Lemma \ref{lem2.1} in order to obtain 
\begin{gather*}
\big\|uv\big\|_{L^{r_1}} \lesssim \big\|u\big\|_{L^{p_1}}
\big\|v\big\|_{L^{q_1}},\ \forall  u \in L^{p_1}, v \in L^{q_1},\\
\big\|uv\big\|_{\dot{H}^1_{r_2}} \lesssim 
\big\|u\big\|_{\dot{H}^1_{p_2}}\big\|v\big\|_{\dot{H}^1_{q_2}}, 
\ \forall  u \in \dot{H}^1_{p_2}, v \in \dot{H}^1_{q_2}.
\end{gather*}
From  Theorem 4.4.1 (page 96) of \cite{J. Bergh: 1976} we get
$$
\big\|uv\big\|_{{\dot{H}^s_r}} \lesssim 
\big\|u\big\|_{\dot{H}^s_p}\big\|v\big\|_{\dot{H}^s_q}.
$$
 \qed
\begin{Lem}\label{lem2.2}
Assume that  
\begin{equation}\label{Them}
q> 1, p>1, 0 \leq \frac{s}{d} < \min\Big\{\frac{1}{p}, \frac{1}{q}\Big\}, 
\ and\ \frac{1}{p} + \frac{1}{q}  <  1 + \frac{s}{d}. 
\end{equation} 
Then we have the inequality 
\begin{gather*}
\big\|uv\big\|_{\dot{H}^s_r} \lesssim 
\big\|u\big\|_{\dot{H}^s_p}\big\|v\big\|_{\dot{H}^s_q},
\ \forall  u \in \dot{H}^s_p,  v \in \dot{H}^s_q,
\end{gather*}
where $\frac{1}{r} = \frac{1}{p} + \frac{1}{q} - \frac{s}{d}$.
\end{Lem}
\textbf{Proof}. Denote by $[s]$ the integer part of $s$ and by $\{s\}$ 
the fraction part of the argument $s$. Using the formula for the 
derivatives of a product of two functions, we have
\begin{gather*}
\big\|uv\big\|_{\dot{H}^s_r} = \big\|\dot{\Lambda}^s(uv)\big\|_{L^r} 
= \big\|\dot{\Lambda}^{\{s\}}(uv)\big\|_{\dot{H}^{[s]}_r} \simeq \\
\sum_{|\alpha| = [s] }\big\|\partial^{\alpha}
\dot{\Lambda}^{\{s\}}(uv)\big\|_{L^r}  = 
\sum_{|\alpha| = [s] }\big\|\dot{\Lambda}^{\{s\}}
\partial^{\alpha}(uv)\big\|_{L^r}\\
= \sum_{|\alpha| = [s] }\big\|\partial^{\alpha}(uv)
\big\|_{\dot{H}^{\{s\}}_r} \lesssim  \sum_{|\gamma| + |\beta| = [s]}
\big\|\partial^{\gamma}u\partial^{\beta}v\big\|_{\dot{H}^{\{s\}}_r}.
\end{gather*}
Set
$$
\frac{1}{\tilde p} = \frac{1}{p} - \frac{s -|\gamma| - \{s\}}{d}, 
\frac{1}{\tilde q} = \frac{1}{q} - \frac{s -|\beta| - \{s\}}{d}.   
$$
Applying Lemma \ref{lem2.1'} and 
the Sobolev inequality in order to obtain 
$$
\big\|\partial^{\gamma}u\partial^{\beta}v\big\|_{\dot{H}^{\{s\}}_r}
 \lesssim \big\|\partial^{\gamma}u\big\|_{\dot{H}^{\{s\}}_{\tilde p}} 
\big\|\partial^{\beta}v\big\|_{\dot{H}^{\{s\}}_{\tilde q}} 
\lesssim \big\|u\big\|_{\dot{H}^{|\gamma| + \{s\}}_{\tilde p}} 
\big\|v\big\|_{\dot{H}^{|\beta| + \{s\}}_{\tilde q}} 
\lesssim \big\|u\big\|_{\dot{H}^s_p}\big\|v\big\|_{\dot{H}^s_q}.
$$
This gives the desired result
$$
\big\|uv\big\|_{\dot{H}^s_r} \lesssim 
\big\|u\big\|_{\dot{H}^s_p}\big\|v\big\|_{\dot{H}^s_q}.
$$
\qed 
\begin{Lem}\label{lem.2.1.4'} 
Let $1 \leq p, q \leq \infty$ and $s \in \mathbb{R}$.\\
{\rm (a)} If $ s < 1$ then the two quantities
\begin{gather}
\ \Big(\int_0^\infty\big(t^{-\frac{s}{2}}\big\|e^{t\Delta}
t^{\frac{1}{2}}\dot{\Lambda}f\big\|_q\big)^p\frac{\dif t}{t} \Big)^{1/p} 
\ and \ \big\|f\big\|_{\dot{B}_{q}^{s, p}}
 \ are \ equivalent. \notag
\end{gather}
{\rm (b)} If $ s < 0$ then the two quantities 
\begin{gather}
\ \Big(\int_0^\infty\big(t^{-\frac{s}{2}}\big\|e^{t\Delta}f\big\|_q\big)^p
\frac{\dif t}{t} \Big)^{1/p} \ and \ \big\|f\big\|_{\dot{B}_{q}^{s, p}}
 \ are \ equivalent, \notag
\end{gather}
where $\dot{B}_{q}^{s, p}$ is the homogeneous Besov space.
\end{Lem}
\textbf{Proof}. See (\cite{S. Friedlander 2004}, Proposition 1, p. 181
 and Proposition 3, p. 182), or see (\cite{P. G. Lemarie-Rieusset 2002}, 
Theorem 5.4, p. 45).\qed\\
The following lemma is a generalization of the above lemma.
\begin{Lem}\label{lem.2.1.4} 
Let $1 \leq p, q \leq \infty,\ \alpha \geq 0$, and $s < \alpha$. 
Then the two quantities
\begin{gather*}
\Big(\int_0^\infty(t^{-\frac{s}{2}}\big\|e^{t\Delta}
t^{\frac{\alpha}{2}} \dot{\Lambda}^{\alpha}f\big\|_{L^q})^p
\frac{{\rm d}t}{t}\Big)^{\frac{1}{p}}\ and 
\ \big\|f\big\|_{\dot{B}_{q}^{s, p}} \ are \ equivalent, 
\end{gather*}
\end{Lem}
\textbf{Proof}. Note that $\dot{\Lambda}^{s_0}$ 
is an isomorphism from $\dot{B}^{s,p}_q$ to $\dot{B}^{s-s_0,p}_q$, 
see \cite{Bjorn Jawerth: 1977}, then we can easily prove the lemma.   \qed
\begin{Lem}\label{lem.3.1} Assume that 
$q > 1, 1 \leq r \leq \infty,\ and\ 0 \leq s < \frac{d}{q}$. The following statement is true: 
If $u_0 \in \dot{H}^s_{L^{q, r}}$ then 
$\me^{t\Delta}u_0 \in L^\infty([0, \infty); \dot{H}^s_{L^{q, r}})$ and 
$\big\|\me^{t\Delta}u_0\big\|_{L^\infty([0, \infty); \dot{H}^s_{L^{q, r}})} 
\leq\big\|u_0\big\|_{\dot{H}^s_{L^{q, r}}}.$
\end{Lem}
\textbf{Proof}. We have
\begin{gather*}
\big\|\me^{t\Delta}u_0\big\|_{\dot{H}^s_{L^{q, r}}} = \big\|\me^{t\Delta}
\dot{\Lambda}^su_0\big\|_{L^{q, r}} 
= \frac{1}{(4\pi t)^{d/2}} \Big\|\int_{\mathbb R^d}\me^{\frac{-|\xi|^2}{4t}}
\dot{\Lambda}^su_0(\ . -\xi)\dif\xi\Big\|_{L^{q, r}} \\ \leq
\frac{1}{(4\pi t)^{d/2}}\int_{\mathbb R^d}\me^{\frac{-|\xi|^2}{4t}}
\big\|\dot{\Lambda}^su_0(\ . -\xi)\big\|_{L^{q, r}}\dif\xi \\
=  \frac{1}{(4\pi t)^{d/2}}\int_{\mathbb R^d}\me^{\frac{-|\xi|^2}{4t}}
\big\|u_0\big\|_{\dot{H}^s_{L^{q, r}}}\dif\xi =
\big\|u_0\big\|_{\dot{H}^s_{L^{q, r}}}.
\end{gather*}
\qed \\
Let us recall following result on solutions of a quadratic
equation in Banach spaces (Theorem 22.4 in 
\cite{P. G. Lemarie-Rieusset 2002}, p. 227).
\begin{Th}\label{th.3.1}
Let $E$ be a Banach space, and $B: E \times E \rightarrow  E$ 
be a continuous bilinear map such that there exists $\eta > 0$ so that
$$
\|B(x, y)\| \leq \eta \|x\| \|y\|,
$$
for all x and y in $E$. Then for any fixed $y \in E$ 
such that $\|y\| \leq \frac{1}{4\eta}$, the equation $x = y - B(x,x)$ 
has a unique solution  $\overline{x} \in E$ satisfying 
$\|\overline{x}\| \leq \frac{1}{2\eta}$.
\end{Th}
\vskip 0.5cm
\centerline{\S3. Main results} 
\vskip 0.5cm
Now, for $T > 0$, we say that $u$ is a mild solution of NSE on $[0, T]$ 
corresponding to a divergence-free initial datum $u_0$ 
when $u$ solves the integral equation
$$
u = e^{t\Delta}u_0 - \int_{0}^{t} e^{(t-\tau) \Delta} \mathbb{P} 
\nabla  .\big(u(\tau,.)\otimes u(\tau,.)\big) \dif\tau.
$$
Above we have used the following notation: 
For a tensor $F = (F_{ij})$ 
we define the vector $\nabla.F$ by 
$(\nabla.F)_i = \sum_{j = 1}^d\partial_jF_{ij}$ 
and for two vectors $u$ and $v$, we define their tensor 
product $(u \otimes v)_{ij} = u_iv_j$. 
The operator $\mathbb{P}$ is the Helmholtz-Leray 
projection onto the divergence-free fields 
\begin{equation}\label{3.1} 
(\mathbb{P}f)_j =  f_j + \sum_{1 \leq k \leq d} R_jR_kf_k, 
\end{equation} 
where $R_j$ is the Riesz transforms defined as 
$$
R_j = \frac{\partial_j}{\sqrt{- \Delta}},\ \ {\rm i.\ e.} \ \  
\widehat{R_jg}(\xi) = \frac{i\xi_j}{|\xi|}\hat{g}(\xi)
$$
with $\hat{}$ denoting the Fourier transform. 
The heat kernel $e^{t\Delta}$ is defined as 
$$
e^{t\Delta}u(x) = ((4\pi t)^{-d/2}e^{-|.|^2/4t}*u)(x).
$$
If $X$ is a normed space and 
$u = (u_1, u_2,...,u_d), u_i \in X, 1 \leq i \leq d$, 
then we write 
$$
u \in X, \|u\|_X = \Big(\sum_{i = 1}^d\|u_i\|_X^2\Big)^{1/2}.
$$
We define the auxiliary space $\mathcal{K}^{s, \tilde q}_{q,r,T}$ 
which is made up by the functions $u(t,x)$ such that 
$$
\big\|u\big\|_{\mathcal{K}^{s, \tilde q}_{q,r,T}}:= 
\underset{0 < t < T}{\rm sup}t^{\frac{\alpha}{2}}
\big\|u(t, .)\big\|_{\dot{H}^s_{L^{\tilde q,r}}} < \infty,
$$
and
\begin{equation}\label{3.1.2}
\underset{t \rightarrow 0}{\rm lim}
t^{\frac{\alpha}{2}}\big\|u(t, .)
\big\|_{\dot{H}^s_{L^{\tilde q,r}}} = 0,
\end{equation}
where $r,q, \tilde q, s$ being fixed constants satisfying 
$$
q, \tilde q  \in (1, + \infty), r \geq 1,  s \geq 0, \frac{s}{d} 
< \frac{1}{\tilde q} \leq \frac{1}{q} \leq \frac{s + 1}{d},
$$ 
and  
$$
\alpha =\alpha(q,\tilde q)  = d\Big(\frac{1}{q} - \frac{1}{\tilde q}\Big).
$$
In the case $\tilde q   = q$, 
it is also convenient to define the space 
$\mathcal{K}^{s,\tilde q}_{q,r,T}$ as 
the natural space $L^\infty([0, T]; \dot{H}^s_{L^{q,r}}(\mathbb{R}^d))$ 
with the additional condition that its elements $u(t,x)$ satisfy 
\begin{equation}\label{3.1.3}
\underset{t \rightarrow 0}{\rm lim}
\big\|u(t,.)\big\|_{\dot{H}^s_{L^{q,r}}} = 0.
\end{equation}
\begin{Not}
The auxiliary space $\mathcal{K}_{\tilde q}: 
= \mathcal{K}^{0, \tilde q}_{d,\tilde q,T}\  (\tilde q \geq d)$ 
was introduced by Weissler and systematically used by 
Kato \cite{T. Kato 1984} and Cannone \cite{M. Cannone 1997}.
\end{Not}
\begin{Lem}\label{lem3.1} Let $1 \leq r \leq \tilde{r}\leq \infty$. 
Then we have the following imbedding maps
$$
\mathcal{K}^{s,\tilde q}_{q,1,T} 
\hookrightarrow \mathcal{K}^{s,\tilde q}_{q,r,T}
 \hookrightarrow \mathcal{K}^{s,\tilde q}_{q,\tilde r,T} \hookrightarrow 
\mathcal{K}^{s,\tilde q}_{q,\infty,T}.
$$
\end{Lem}
\textbf{Proof}. It is easily deduced from Lemma \ref{lem2.1.1} (a) 
and the definition \linebreak of  $\mathcal{K}^{s,\tilde q}_{q,r,T}$. \qed

\begin{Lem}\label{lem.3.5} If \ $u_0 
\in \dot{H}^s_{L^{q,r}}(\mathbb{R}^d)$ with 
$q >1, r \geq 1, s \geq 0$, and $ \frac{s}{d} 
<  \frac{1}{q} \leq \frac{s + 1}{d}$ then for all $\tilde q $ satisfying
$$
\frac{s}{d} < \frac{1}{\tilde q} < \frac{1}{q},
$$
we have
$$
e^{t\Delta}u_0 \in \mathcal{K}^{s, \tilde q}_{q,1,\infty}, 
$$
and the following imbedding map
\begin{equation}\label{lem.3.3.1.a}
\dot{H}^s_{L^{q,r}}(\mathbb{R}^d) \hookrightarrow   \dot{B}^{s- (\frac{d}{q} 
- \frac{d}{\tilde q}), \infty}_{\tilde q}(\mathbb{R}^d).
\end{equation}
\end{Lem}
\textbf{Proof}. Before proving this lemma, we need 
to prove the following lemma.
\begin{Lem}\label{lem.3.4.1} Suppose that 
$u_0 \in  L^{q, r}(\mathbb{R}^d)$ with $1 \leq q \leq \infty$ 
and $1 \leq r < \infty$. Then 
$\underset{n \rightarrow \infty}{\rm lim}
\big\|\mathcal{X}_nu_0\big\|_{L^{q, r}} = 0$,
where  $n \in \mathbb{N}, \mathcal{X}_n(x) = 0$ for $x \in \{x : \ |x| < n\} 
\cap \{x : \big|u_0(x)\big| < n\}$ 
and $\mathcal{X}_n(x) = 1$ otherwise.
\end{Lem}
\textbf{Proof}.  With $\delta > 0$ being fixed, we have
\begin{equation}\label{lem.3.3.1}
\big\{x: |\mathcal{X}_nu_0(x)| > \delta\big\} \supseteq 
\big\{x: |\mathcal{X}_{n+1}u_0(x)| > \delta\big\},
\end{equation}
and
\begin{equation}\label{lem.3.3.2}
\underset{n = 0}{\overset{\infty }{\cap}}
\{x: |\mathcal{X}_nu_0(x)| > \delta\} = \{x: |u_0(x)| = + \infty\}.
\end{equation}
We prove that
\begin{equation}\label{lem.3.3.2.a}
\mathcal{M}^d (\{x: |u_0(x)| = + \infty\}) = 0,
\end{equation}
with $\mathcal{M}^d$ being the  Lebesgue measure in $\mathbb R^d$, assuming on the contrary
$$
\mathcal{M}^d (\{x : |u_0(x)| = + \infty\}) > 0.
$$
We have $u_0^*(t): = \inf \big\{\tau : \mathcal{M}^d
\big(\{x: |u_0(x)| > \tau\}\big) \leq t \big\} = +\infty$ for all $t$ such that $0 < t < \mathcal{M}^d (\{x : |u_0(x)| = + \infty\})$ and then $\big\|u_0\big\|_{L^{q,r}} = + \infty$, a contradiction.\\
Note that 
$$
\mathcal{M}^d\big(\{x: |\mathcal{X}_0u_0(x)| > \delta\}\big) 
= \mathcal{M}^d\big(\{x: |u_0(x)| > \delta\}\big).
$$
We prove that
\begin{equation}\label{lem.3.3.3}
\mathcal{M}^d\big(\{x: |u_0(x)| > \delta\}\big) < \infty,
\end{equation}
assuming on the contrary
$$
\mathcal{M}^d\big(\{x: |u_0(x)| > \delta\}\big) = \infty.
$$
We have $u_0^*(t) \geq \delta\ \text{for all}\ t > 0$, from the definition 
of the Lorentz space, we get 
$$
\big\|u_0\big\|_{L^{q,r}} = \Big(\int_0 ^\infty (t^{\frac{1}{q}}u_0^*(t))^r
\frac{\dif t}{t}\Big)^{\frac{1}{r}} \geq \Big(\int_0 ^\infty (t^{\frac{1}{q}}\delta)^r
\frac{\dif t}{t}\Big)^{\frac{1}{r}} =\delta \Big(\int_0 ^\infty t^{\frac{r}{q} - 1}
\dif t\Big)^{\frac{1}{r}} = \infty,
$$
a contradiction.\\
From \eqref{lem.3.3.1}, \eqref{lem.3.3.2}, \eqref{lem.3.3.2.a},
and \eqref{lem.3.3.3}, we infer that 
\begin{equation}\label{lem.3.3.4}
\underset{n \rightarrow \infty}{\rm lim}\mathcal{M}^d
\big(\{x: |\mathcal{X}_nu_0(x)| > \delta\}\big) = \mathcal{M}^d (\{x: |u_0(x)| = + \infty\}) = 0.
\end{equation}
Set
$$
u^*_n(t) = \inf \big\{\tau : \mathcal{M}^d
\big(\{x: |\mathcal{X}_nu_0(x)| > \tau\}\big) \leq t \big\}.
$$
We have
\begin{equation}\label{lem.3.3.5}
u^*_n(t) \geq u^*_{n + 1}(t).
\end{equation}
Fixed $t > 0$. For any $\epsilon > 0$, from \eqref{lem.3.3.4} it follows that 
there exists a number $n_0 = n_0(t, \epsilon)$ large enough such that 
$$
\mathcal{M}^d\big(\{x: |\mathcal{X}_nu_0(x)| 
> \epsilon\}\big) \leq t , \forall n \geq n_0.
$$
From this we deduce that
$$
u^*_n(t) \leq \epsilon , \forall n \geq n_0,
$$
therefore
\begin{equation}\label{lem.3.3.6}
\lim_{n\rightarrow \infty} u^*_n(t) = 0.
\end{equation}
From \eqref{lem.3.3.5} and \eqref{lem.3.3.6}, we apply 
Lebesgue's monotone convergence theorem to get
$$
\underset{n \rightarrow \infty}{\rm lim}\big\|\mathcal{X}_nu_0\big\|_{L^{q, r}} =
 \underset{n \rightarrow \infty}{\rm lim}\Big(\int_0 ^\infty 
(t^{\frac{1}{q}}u^*_n(t))^r\frac{\dif t}{t}\Big)^{\frac{1}{r}} = 0. \qed
$$
Now we return to prove Lemma \ref{lem.3.5}. We prove that
\begin{equation}\label{lem.3.3.1.b}
\underset{0 < t < \infty}{\rm sup}t^{\frac{\alpha}{2}}
\big\|e^{t\Delta}u_0\big\|_{\dot{H}^s_{L^{\tilde q,1}}} 
\lesssim \big\|u_0\big\|_{\dot{H}^s_{L^{q,r}}}.
\end{equation}
Set
$$
\frac{1}{h} = 1 +  \frac{1}{\tilde q} - \frac{1}{q}.
$$
Applying Proposition 2.4 $(c)$ in 
(\cite{P. G. Lemarie-Rieusset 2002}, pp. 20) 
for convolution in the Lorentz spaces, we have
\begin{gather*}
\big\|e^{t\Delta}u_0\big\|_{\dot{H}^s_{L^{\tilde q,1}}} 
= \Big\|e^{t\Delta}\dot{\Lambda}^su_0\Big\|_{L^{\tilde q,1}} =
\frac{1}{(4\pi t)^{d/2}} \Big\|e^{-\frac{| . |^2}
{4t}}*\dot{\Lambda}^su_0\Big\|_{L^{\tilde q,1}} \notag \lesssim \\
\frac{1}{t^{d/2}}\big\|e^{-\frac{| . |^2}{4t}}\big\|_{L^{h,1}}
\big\|\dot{\Lambda}^su_0\big\|_{L^{q,\infty}} = 
t^{-\frac{\alpha}{2}}\big\|e^{-\frac{| . |^2}{4}}\big\|_{L^{h,1}}
\big\|u_0\big\|_{\dot{H}^s_{L^{q,\infty}}}  
\lesssim t^{-\frac{\alpha}{2}}\big\|u_0\big\|_{\dot{H}^s_{L^{q,r}}}. 
\end{gather*}
We claim now that
$$
\underset{t \rightarrow 0}{\rm lim}t^{\frac{\alpha}{2}}
\big\|e^{t\Delta}u_0\big\|_{\dot{H}^s_{L^{\tilde q,1}}} = 0.
$$
From Lemma \ref{lem.3.4.1}, we have
\begin{equation}\label{3.8}
\underset{n \rightarrow \infty}{\rm lim}
\Big\|\mathcal{X}_{n,s}\dot{\Lambda}^su_0\Big\|_{L^{q,r}} = 0,
\end{equation}
where $\mathcal{X}_{n,s}(x) = 0$ for $x \in \{x : \ |x| < n\} 
\cap \{x : \big|\dot{\Lambda}^su_0(x)\big| < n\}$ 
and $\mathcal{X}_{n,s}(x) = 1$ otherwise. 
We have
\begin{gather}
t^{\frac{\alpha}{2}}\big\|e^{t\Delta}u_0\big\|_{\dot{H}^s_{L^{\tilde q,1}}} \leq 
\frac{t^{\frac{\alpha}{2} - \frac{d}{2}}}{(4\pi)^{d/2}}
\Big\|e^{-\frac{| . |^2}{4t}}*(\mathcal{X}_{n,s}\dot{\Lambda}^su_0)
\Big\|_{L^{\tilde q,1}} + \notag \\ 
\frac{t^{\frac{\alpha}{2} - \frac{d}{2}}}{(4\pi)^{d/2}}
\Big\|e^{-\frac{| . |^2}{4t}}*((1 - \mathcal{X}_{n,s})\dot{\Lambda}^su_0)
\Big\|_{L^{\tilde q,1}}.\label{3.7}
\end{gather}
For any $\epsilon > 0$, applying Proposition 2.4 $(c)$ 
in (\cite{P. G. Lemarie-Rieusset 2002}, pp. 20) 
and note that \eqref{3.8}, we have
\begin{gather}
\frac{t^{\frac{\alpha}{2} - \frac{d}{2}}}{(4\pi)^{d/2}}
\Big\|e^{-\frac{| . |^2}{4t}}*(\mathcal{X}_{n,s}
\dot{\Lambda}^su_0)\Big\|_{L^{\tilde q,1}} \notag \\ 
\leq C_1\big\|e^{-\frac{| . |^2}{4}}\big\|_{L^{h,1}}
\Big\|\mathcal{X}_{n,s}\dot{\Lambda}^su_0\Big\|_{L^{q,\infty}} 
\leq C_2\Big\|\mathcal{X}_{n,s}\dot{\Lambda}^su_0
\Big\|_{L^{q,r}} < \frac{\epsilon}{2}, \label{3.9}
\end{gather}
for large enough $n$. Fixed one of such $n$, applying  
Proposition 2.4 $(a)$ in (\cite{P. G. Lemarie-Rieusset 2002}, 
pp. 20), we conclude that 
\begin{gather}
\frac{t^{\frac{\alpha}{2} - \frac{d}{2}}}{(4\pi)^{d/2}}
\Big\|e^{-\frac{| . |^2}{4t}}*((1 - \mathcal{X}_{n,s})
\dot{\Lambda}^su_0)\Big\|_{L^{\tilde q,1}} \notag \\
\leq C_3t^{\frac{\alpha}{2} - \frac{d}{2}}
\big\|e^{-\frac{| . |^2}{4t}}\big\|_{L^1}
\Big\|(1 - \mathcal{X}_{n,s})\dot{\Lambda}^su_0
\Big\|_{L^{\tilde q,1}}  \notag \\
 \leq C_4t^{\frac{\alpha}{2}}\big\|e^{-\frac{| . |^2}{4}}
\big\|_{L^1}\big\|n(1 - \mathcal{X}_{n,s})\big\|_{L^{\tilde q,1}} = \notag \\
C_5nt^{\frac{\alpha}{2}}\big\|(1 - \mathcal{X}_{n,s})\big\|_{L^{\tilde q,1}} 
= C_6(n)t^{\frac{\alpha}{2}} < \frac{\epsilon}{2},  \label{3.10}
\end{gather}
for small enough $t > 0$.  From the estimates \eqref{3.7}, \eqref{3.9}, 
and \eqref{3.10} it follows that 
$$
t^{\frac{\alpha}{2}}\big\|e^{t\Delta}u_0
\big\|_{\dot{H}^s_{L^{\tilde q,1}}}
 \leq C_2\Big\|\mathcal{X}_{n,s}\dot{\Lambda}^su_0\Big\|_{L^{q,r}} + 
C_6(n)t^{\frac{\alpha}{2}} < \epsilon. 
$$
Finally, the embedding \eqref{lem.3.3.1.a}  is derived from the inequality \eqref{lem.3.3.1.b}, Lemma \ref{lem2.1.1}, and Lemma \ref{lem.2.1.4}.
\begin{Not}
In the case  $s = 0$ and $q=r =d$, Lemma \ref{lem.3.3} 
is a generalization of Lemma $9$ in (\cite{M. Cannone. 2004}, p. 196). 
\end{Not}
In the following lemmas a particular attention 
will be devoted to study of the bilinear operator 
$B(u, v)(t)$ defined by 
$$
B(u, v)(t) = \int_{0}^{t} e^{(t-\tau ) \Delta} \mathbb{P} 
\nabla.\big(u(\tau)\otimes v(\tau)\big) \dif\tau.
$$
\begin{Lem}\label{lem.3.3}
Let $s, q \in \mathbb R$ be such that
\begin{equation}\label{3.11.a}
s \geq 0, q > 1 ,\ and\ \frac{s}{d} <  \frac{1}{q} 
\leq \frac{s + 1}{d}. 
\end{equation}
Then for all $\tilde q$ satisfying
\begin{equation}\label{3.11}
\frac{s}{d} < \frac{1}{\tilde q} < {\rm min} 
\Big\{\frac{1}{2} +\frac{s}{2d}, \frac{1}{q}\Big\},
\end{equation}
the bilinear operator $B(u, v)(t) $ is continuous from 
$\mathcal{K}_{q,\tilde q,T}^{s, \tilde q} 
\times \mathcal{K}_{q,\tilde q,T}^{s, \tilde q}$ into 
$\mathcal{K}_{q,1,T}^{s, \tilde q}$ and the following inequality holds 
\begin{equation}\label{3.12}
\big\|B(u, v)\big\|_{\mathcal{K}_{q,1,T}^{s, \tilde q}} \leq
 C.T^{\frac{1}{2}(1 + s -\frac{d}{q})}
\big\|u\big\|_{\mathcal{K}_{q,\tilde q,T}^{s, \tilde q}}
\big\|v\big\|_{\mathcal{K}_{q,\tilde q,T}^{s, \tilde q}},
\end{equation}
where C is a positive constant independent of T.
\end{Lem}
\textbf{Proof}. We have 
\begin{gather}
\big\|B(u,v)(t)\big\|_{\dot{H}^s_{L^{\tilde q,1}}}
 \leq \int_{0}^{t} \Big\|e^{(t- \tau) \Delta}
 \mathbb{P} \nabla .\big(u(\tau,.)\otimes v(\tau,.)\big)
\Big\|_{\dot{H}^s_{L^{\tilde q,1}}} d \tau  = \notag \\
 \int_{0}^{t} \Big\|e^{(t- \tau) \Delta} \mathbb{P} \nabla .\dot{\Lambda}^s
\big(u(\tau,.)\otimes v(\tau,.)\big)\Big\|_{L^{\tilde q,1}}d \tau.\label{3.13}
\end{gather}
From the properties of the Fourier transform
\begin{gather}
\Big(e^{(t - \tau)\Delta}\mathbb{P}\nabla .\dot{\Lambda}^s
\big(u(\tau,.) \otimes v(\tau,.)\big)\Big)_j^\wedge (\xi) = \notag \\  
e^{-(t - \tau)|\xi|^2}\sum_{l, k = 1}^d \Big(\delta_{jk} - 
\frac{\xi_j\xi_k}{|\xi|^2}\Big)(i\xi_l)\Big(\dot{\Lambda}^s
\big(u_l(\tau,.)v_k(\tau,.)\big)\Big)^\wedge(\xi), \notag
\end{gather}
and then
\begin{gather}
\Big(e^{(t - \tau)\Delta}\mathbb{P}\nabla .\dot{\Lambda}^s
\big(u(\tau,.) \otimes v(\tau,.)\big)\Big)_j = \notag \\
 \frac{1}{(t - \tau)^{\frac{d + 1}{2}}}\sum_{l, k = 1}^d K_{l, k, j}
\Big(\frac{.}{\sqrt{t - \tau}}\Big)*\Big(\dot{\Lambda}^s
\big(u_l(\tau,.)v_k(\tau,.)\big)\Big), \label{3.14}
\end{gather}
where 
$$
\widehat{K_{l, k, j}}(\xi)= \frac{1}{(2\pi)^{d/2}}.e^{-|\xi|^2}
\Big(\delta_{jk} - \frac{\xi_j\xi_k}{|\xi|^2}\Big)(i\xi_l).
$$
Applying Proposition 11.1 (\cite{P. G. Lemarie-Rieusset 2002}, p. 107)  
with $|\alpha| = 1$ we see that the tensor $K(x) = \{K_{l, k, j}(x)\}$ satisfies 
\begin{equation}\label{3.15}
|K(x)| \lesssim \frac{1}{(1 + |x|)^{d + 1}}.
\end{equation}
So, we can rewrite the equality \eqref{3.14}  in the tensor form 
$$
e^{(t - \tau)\Delta}\mathbb{P}\nabla .\dot{\Lambda}^s\big(u(\tau,.) 
\otimes v(\tau,.)\big)  = 
$$
\begin{equation}\label{3.17.a}
\frac{1}{(t - \tau)^{\frac{d + 1}{2}}}K\Big(\frac{.}{\sqrt{t - \tau}}\Big)*
\Big(\dot{\Lambda}^s\big(u(\tau,.) \otimes v(\tau,.)\big)\Big).
\end{equation}
Set
\begin{equation}\label{3.17}
\frac{1}{r} = \frac{2}{\tilde{q}} - \frac{s}{d},\ \frac{1}{h} 
= \frac{s}{d} - \frac{1}{\tilde{q}} + 1.
\end{equation}
From the inequalities \eqref{3.11.a} and \eqref{3.11}, 
we can check that the following conditions are satisfied
$$
1 < h, r < \infty \ \text{and}\  \frac{1}{\tilde q} 
+ 1 = \frac{1}{h} + \frac{1}{r}.
$$
Applying Proposition 2.4 $(c)$ in (\cite{P. G. Lemarie-Rieusset 2002}, 
pp. 20) for convolution in the Lorentz spaces, we have
\begin{gather}
\Big\|e^{(t - \tau)\Delta}\mathbb{P}\nabla .\dot{\Lambda}^s\big(u(\tau,.) 
\otimes v(\tau,.)\big)\Big\|_{L^{\tilde q,1}} \lesssim \notag \\
\frac{1}{(t - \tau)^{\frac{d + 1}{2}}}
\Big\|K\Big(\frac{.}{\sqrt{t - \tau}}\Big)\Big\|_{L^{h,1}}
\Big\|\dot{\Lambda}^s\big(u(\tau,.) \otimes v(\tau,.)\big)
\Big\|_{L^{r,\infty}}. \label{3.16}
\end{gather}
Applying Lemma \ref{lem2.2} we obtain
\begin{gather}
\Big\|\dot{\Lambda}^s\big(u(\tau,.) \otimes v(\tau,.)\big)
\Big\|_{L^{r,\infty}} \leq \Big\|\dot{\Lambda}^s\big(u(\tau,.) 
\otimes v(\tau,.)\big)\Big\|_{L^r}
= \big\|u(\tau,.)\otimes v(\tau,.)\big\|_{\dot{H}^s_r} \notag \\
\lesssim \big\|u(\tau,.)\big\|_{\dot{H}^s_{\tilde q}}\big\|v(\tau,.)
\big\|_{\dot{H}^s_{\tilde q}}.
\label{3.18}
\end{gather}
Fom the inequalities \eqref{3.15} and \eqref{3.17} we infer that 
\begin{gather}
\Big\|K\Big(\frac{.}{\sqrt{t - \tau}}\Big)\Big\|_{L^{h,1}} 
= (t - \tau)^{\frac{d}{2h}}\big\|K\big\|_{L^{h,1}} \simeq (t - \tau)^{\frac{s}{2} 
- \frac{d}{2{\tilde q}} + \frac{d}{2}}.
\label{3.19}
\end{gather}
From the inequalities \eqref{3.16}, \eqref{3.18}, and \eqref{3.19} 
we deduce that 
\begin{gather}
\Big\|e^{(t - \tau)\Delta}\mathbb{P}\nabla .\dot{\Lambda}^s\big(u(\tau,.) 
\otimes v(\tau,.)\big)\Big\|_{L^{\tilde q,1}} \lesssim  
\notag \\  (t - \tau)^{\frac{s}{2} - \frac{d}{2{\tilde q}} 
- \frac{1}{2}}\big\|u(\tau,.)\big\|_{\dot{H}^s_{\tilde q}}
\big\|v(\tau,.)\big\|_{\dot{H}^s_{\tilde q}}. \label{3.20}
\end{gather}
From the estimates \eqref{3.13} and \eqref{3.20}, and 
note that from the inequalities \eqref{3.11.a} and \eqref{3.11}, 
we can check that $\frac{s}{2} - \frac{d}{2{\tilde q}}
 - \frac{1}{2} > -1$ and $\alpha = d(\frac{1}{q}
 - \frac{1}{\tilde q}) < 1$, this gives the desired result
\begin{gather}
\big\|B(u, v)(t)\big\|_{\dot{H}^s_{L^{\tilde q,1}}} \lesssim 
\int_0^t (t - \tau)^{\frac{s}{2} - \frac{d}{2{\tilde q}}
- \frac{1}{2}}\big\|u(\tau,.)\big\|_{\dot{H}^s_{\tilde q}}.
\big\|v(\tau,.)\big\|_{\dot{H}^s_{\tilde q}}{\rm d} \tau \lesssim \notag \\
\int_0^t (t - \tau)^{\frac{s}{2} - \frac{d}{2{\tilde q}}- \frac{1}{2}}\tau^{-\alpha}
\underset{0 < \eta < t}{\rm sup}\eta^{\frac{\alpha}{2}}
\big\|u(\eta, .)\big\|_{\dot{H}^s_{\tilde q}}.\underset{0 < \eta < t}{\rm sup}
\eta^{\frac{\alpha}{2}}\big\|v(\eta, .)\big\|_{\dot{H}^s_{\tilde q}}\dif\tau = \notag \\
\underset{0 < \eta < t}{\rm sup}\eta^{\frac{\alpha}{2}}\big\|u(\eta, .)
\big\|_{\dot{H}^s_{\tilde q}}.\underset{0 < \eta < t} {\rm sup}
\eta^{\frac{\alpha}{2}}\big\|v(\eta, .)\big\|_{\dot{H}^s_{\tilde q}}
\int_0^t (t - \tau)^{\frac{s}{2} - \frac{d}{2{\tilde q}}
- \frac{1}{2}}\tau^{-\alpha}\dif\tau \simeq \notag \\
t^{-\frac{\alpha}{2}}t^{\frac{1}{2}(1 + s -\frac{d}{q})}
\underset{0 < \eta < t}{\rm sup}\eta^{\frac{\alpha}{2}}
\big\|u(\eta, .)\big\|_{\dot{H}^s_{L^{\tilde q,\tilde q}}}.
\underset{0 < \eta < t}{\rm sup}
\eta^{\frac{\alpha}{2}}\big\|v(\eta, .)
\big\|_{\dot{H}^s_{L^{\tilde q,\tilde q}}}. \label{3.21}
\end{gather}
Let us now check the validity of the condition \eqref{3.1.2} 
for the bilinear term $B(u,v)(t)$.
Indeed, we have 
$$
\underset{t \rightarrow 0}{\rm lim}t^{\frac{\alpha}{2}}
\big\|B(u,v)(t)\big\|_{\dot{H}^s_{L^{\tilde q,1}}} = 0,
$$
whenever
$$
\underset{t \rightarrow 0}{\rm lim}t^{\frac{\alpha}{2}}
\big\|u(t,.)\big\|_{\dot{H}^s_{\tilde q}} 
= \underset{t \rightarrow 0}{\rm lim}t^{\frac{\alpha}{2}}
\big\|v(t,.)\big\|_{\dot{H}^s_{\tilde q}} = 0.
$$
The estimate  \eqref{3.12} is now deduced from 
the inequality \eqref{3.21}. \qed 
\begin{Not}
In the case $s = 0$ and $q=d$, Lemma \ref{lem.3.5} is 
a generalization of Lemma $10$ in (\cite{M. Cannone. 2004}, p. 196). 
\end{Not}
\begin{Lem}\label{lem.3.4}
Let $s, q \in \mathbb R$ be such that
\begin{equation}\label{3.24}
s \geq 0, q > 1,\ and\ \frac{s}{d} <  \frac{1}{q} 
\leq \frac{s + 1}{d}. 
\end{equation}
Then for all $\tilde q$ satisfying
\begin{equation}\label{3.24.a}
\frac{1}{2}\Big({\frac{1}{q}+\frac{s}{d}}\Big) 
< \frac{1}{\tilde q} < {\rm min} 
\Big\{\frac{1}{2} +\frac{s}{2d}, \frac{1}{q}\Big\}, 
\end{equation}
the bilinear operator $B(u, v)(t) $ is 
continuous from $\mathcal{K}_{q,\tilde q,T}^{s,\tilde q} 
\times \mathcal{K}_{q,\tilde q,T}^{s,\tilde q}$ into 
$\mathcal{K}_{q,1,T}^{s,q}$ and the following inequality holds 
\begin{equation}\label{3.25}
\big\|B(u, v)\big\|_{\mathcal{K}_{q,1,T}^{s,q}} \leq
  C.T^{\frac{1}{2}(1 + s -\frac{d}{q})}
\big\|u\big\|_{\mathcal{K}_{q,\tilde q,T}^{s,\tilde q}}
\big\|v\big\|_{\mathcal{K}_{q,\tilde q,T}^{s,\tilde q}},
\end{equation}
where C is a positive constant independent of T.
\end{Lem}
\textbf{Proof}. Set
\begin{equation}\label{3.27}
\frac{1}{r} = \frac{2}{\tilde q} - \frac{s}{d},\ \frac{1}{h} 
= 1 +  \frac{1}{q} - \frac{2}{\tilde q} + \frac{s}{d}.
\end{equation}
From the inequalities \eqref{3.24} and \eqref{3.24.a}, 
we can check that $h$ and $r$ satisfy 
$$
1 < h, r < \infty\ \text{and}\ \frac{1}{q} + 1 = \frac{1}{h} + \frac{1}{r}. 
$$ 
From the equality \eqref{3.17.a}, 
applying Proposition 2.4 $(c)$ in (\cite{P. G. Lemarie-Rieusset 2002}, 
pp. 20), we obtain 
\begin{gather}
\Big\|e^{(t - \tau)\Delta}\mathbb{P}\nabla .\dot{\Lambda}^s\big(u(\tau,.) 
\otimes v(\tau,.)\big)\Big\|_{L^{q,1}} \lesssim \notag \\
\frac{1}{(t - \tau)^{\frac{d + 1}{2}}}
\Big\|K\Big(\frac{.}{\sqrt{t - \tau}}\Big)\Big\|_{L^{h,1}}
\Big\|\dot{\Lambda}^s\Big(u(\tau,.) 
\otimes v(\tau,.)\Big)\Big\|_{L^{r,\infty}}.\label{3.26}
\end{gather}
Applying Lemma \ref{lem2.2}, we have
\begin{gather}
\Big\|\dot{\Lambda}^s\big(u(\tau,.) \otimes v(\tau,.)\big)
\Big\|_{L^{r,\infty}} \leq \Big\|\dot{\Lambda}^s\big(u(\tau,.) 
\otimes v(\tau,.)\big)\Big\|_{L^r} \notag \\
\lesssim \big\|u(\tau,.)\big\|_{\dot{H}^s_{\tilde q}}
\big\|v(\tau,.)\big\|_{\dot{H}^s_{\tilde q}}.
\label{3.28}
\end{gather}
From the inequalities \eqref{3.15} and \eqref{3.27} it follows that
\begin{gather}
\Big\|K\Big(\frac{.}{\sqrt{t - \tau}}\Big)\Big\|_{L^{h,1}} 
=  (t - \tau)^{\frac{d}{2h}}\big\|K\big\|_{L^{h,1}} \simeq (t - \tau)^{\frac{d}{2}
+ \frac{d}{2q} - \frac{d}{\tilde q} + \frac{s}{2}}.
\label{3.29}
\end{gather}
From the estimates \eqref{3.26}, \eqref{3.28}, \eqref{3.29} we deduce that 
\begin{gather*}
\Big\|e^{(t - \tau)\Delta}\mathbb{P}\nabla .\big(u(\tau,.) 
\otimes v(\tau,.)\big)\Big\|_{\dot{H}^s_{L^{q,1}}} 
\lesssim (t - \tau)^{\frac{d}{2q} - \frac{d}{\tilde q}
 + \frac{s}{2} - \frac{1}{2}}
\big\|u(\tau, .)\big\|_{\dot{H}^s_{\tilde q}}
\big\|v(\tau, .)\big\|_{\dot{H}^s_{\tilde q}} \notag \\
=  (t - \tau)^{\alpha + \frac{s}{2} - \frac{d}{2q} 
- \frac{1}{2}}\big\|u(\tau, .)\big\|_{\dot{H}^s_{\tilde q}}
\big\|v(\tau, .)\big\|_{\dot{H}^s_{\tilde q}}.  
\end{gather*}
From the inequalities \eqref{3.24} and \eqref{3.24.a}, 
we can check that $\alpha + \frac{s}{2} 
- \frac{d}{2q} - \frac{1}{2} > -1$ and 
$\alpha = d(\frac{1}{q} - \frac{1}{\tilde q}) < 1$, 
this gives the desired result
\begin{gather}
\big\|B(u, v)(t)\big\|_{\dot{H}^s_{L^{q,1}}}
\lesssim  \int_0^t (t - \tau)^{\alpha + \frac{s}{2} 
- \frac{d}{2q} - \frac{1}{2}}\big\|u(\tau, .)\big\|_{\dot{H}^s_{\tilde q}}
\big\|v(\tau, .)\big\|_{\dot{H}^s_{\tilde q}}\dif\tau \lesssim \notag \\
\int_0^t (t - \tau)^{\alpha + \frac{s}{2} - \frac{d}{2q} 
- \frac{1}{2}}\tau^{-\alpha}\underset{0 < \eta < t}
{\rm sup}\eta^{\frac{\alpha}{2}}\big\|u(\eta, .)
\big\|_{\dot{H}^s_{\tilde q}}.\underset{0 < \eta < t}{\rm sup}
\eta^{\frac{\alpha}{2}}\big\|v(\eta, .)\big\|_{\dot{H}^s_{\tilde q}}
\dif\tau = \notag \\
\underset{0 < \eta < t}
{\rm sup}\eta^{\frac{\alpha}{2}}\big\|u(\eta, .)
\big\|_{\dot{H}^s_{\tilde q}}.\underset{0 < \eta < t}{\rm sup}
\eta^{\frac{\alpha}{2}}\big\|v(\eta, .)
\big\|_{\dot{H}^s_{\tilde q}}\int_0^t (t - \tau)^{\alpha 
+ \frac{s}{2} - \frac{d}{2q} - \frac{1}{2}}\tau^{-\alpha}\dif\tau
 \simeq \notag \\
t^{\frac{1}{2}{(1 + s - \frac{d}{q}})}\underset{0 < \eta < t}
{\rm sup}\eta^{\frac{\alpha}{2}}
\big\|u(\eta, .)\big\|_{\dot{H}^s_{L^{\tilde q,\tilde q}}}.
\underset{0 < \eta < t}{\rm sup}\eta^{\frac{\alpha}{2}}
\big\|v(\eta, .)\big\|_{\dot{H}^s_{L^{\tilde q,\tilde q}}}. \label{3.30}
\end{gather}
Let us now check the validity of the condition \eqref{3.1.3} 
for the bilinear term $B(u,v)(t)$. 
Indeed, we have 
$$
\underset{t \rightarrow 0}{\rm lim}\big\|B(u,v)(t)
\big\|_{\dot{H}^s_{L^{q,1}}} = 0
$$
whenever
$$
\underset{t \rightarrow 0}{\rm lim}t^{\frac{\alpha}{2}}
\big\|u(t,.)\big\|_{\dot{H}^s_{\tilde q}} 
= \underset{t \rightarrow 0}{\rm lim}t^{\frac{\alpha}{2}}
\big\|v(t,.)\big\|_{\dot{H}^s_{\tilde q}} = 0.
$$
The estimate  \eqref{3.25} is now deduced from 
the inequality \eqref{3.30}. \qed \\
Combining Theorem \ref{th.3.1} with 
Lemmas \ref{lem.3.1}, \ref{lem.3.5}, \ref{lem.3.3}, \ref{lem.3.4},
we obtain the following existence result.
\begin{Th}\label{th.3.2}
Let $s, q,\  and\ r \in \mathbb R$ be such that
\begin{equation}\label{3.33}
s \geq 0, q > 1, r \geq 1, \ and\ \frac{s}{d} <  \frac{1}{q} 
\leq \frac{s + 1}{d}. 
\end{equation}
{\rm (a)} For all $\tilde q$ satisfying
\begin{equation}\label{3.33.a}
\frac{1}{2}\Big({\frac{1}{q}+\frac{s}{d}}\Big) 
< \frac{1}{\tilde q} < {\rm min} \Big\{\frac{1}{2} 
+\frac{s}{2d}, \frac{1}{q}\Big\}, 
\end{equation}
there exists a positive constant $\delta_{s,q,\tilde q,d}$ 
such that for all $T > 0$ and for all 
$u_0 \in \dot{H}^s_{L^{q,r}}(\mathbb{R}^d)$ 
with ${\rm div}(u_0) = 0$ satisfying
\begin{equation}\label{3.34}
T^{\frac{1}{2}(1 + s -\frac{d}{q})}
\underset{0 < t < T}{\rm sup}t^{\frac{d}{2}
({\frac{1}{q} - \frac{1}{\tilde q}})}
\big\|e^{t\Delta}u_0\big\|_{\dot{H}^s_{\tilde q}}
 \leq \delta_{s,q,\tilde q,d},
\end{equation}
NSE has a unique mild solution $u \in \mathcal{K}_{q,1,T}^{s, \tilde q} 
\cap L^\infty\big([0, T]; \dot{H}^s_{L^{q,r}}\big)$. 
In particular, for arbitrary $u_0\in \dot{H}^s_{L^{q,r}}$ with ${\rm div}(u_0) = 0$, 
there exists $T(u_0)$ small enough such that the inequality \eqref{3.34} holds. \\
{\rm (b)} If  $1 < q \leq d ,\ and\ s = \frac{d}{q} - 1$ 
then for any $\tilde q$ be such that 
$$
\frac{1}{q} - \frac{1}{2d} < \frac{1}{\tilde q} 
< {\rm min}\Big\{\frac{1}{2} + \frac{1}{2q} 
- \frac{1}{2d}, \frac{1}{q}\Big\},
$$
there exists a positive constant 
$\sigma_{q,\tilde q, d}$ such that if  
$\big\|u_0\big\|_{\dot{B}^{\frac{d}{\tilde q} - 1, \infty}_{\tilde q}} 
\leq \sigma_{q,\tilde q,d}$ \linebreak $and\ T = \infty$ 
then the inequality \eqref{3.34} holds.
\end{Th}
\textbf{Proof}. From Lemmas \ref{lem.3.3} 
and \ref{lem3.1} , the bilinear operator $B(u,v)(t)$ is continuous 
from $\mathcal{K}_{q,\tilde q,T}^{s, \tilde q} 
\times \mathcal{K}_{q,\tilde q,T}^{s, \tilde q}$ into 
$\mathcal{K}_{q,\tilde q,T}^{s, \tilde q}$ and we have the inequality
$$
\big\|B(u, v)\big\|_{\mathcal{K}_{q,\tilde q,T}^{s, \tilde q}}
 \leq \big\|B(u, v)\big\|_{\mathcal{K}_{q,1,T}^{s, \tilde q}} \leq
 C_{s,q,\tilde q,d}T^{\frac{1}{2}(1 + s -\frac{d}{q})}
\big\|u\big\|_{\mathcal{K}_{q,\tilde q,T}^{s, \tilde q}}
\big\|v\big\|_{\mathcal{K}_{q,\tilde q,T}^{s, \tilde q}},
$$
where $C_{s,q,\tilde q,d}$ is a positive constant independent of  $T$. 
From Theorem \ref{th.3.1} and the above inequality, 
we deduce following: for any $u_0 \in \dot{H}^s_{L^{q,r}}(\mathbb{R}^d)$ 
such that
$$
{\rm div}(u_0)=0,\ T^{\frac{1}{2}(1 + s -\frac{d}{q})}
\underset{0 < t < T}{\rm sup}t^{\frac{d}{2}({\frac{1}{q} - \frac{1}{\tilde q}})}
\big\|e^{t\Delta}u_0\big\|_{\dot{H}^s_{\tilde q}} \leq  \frac{1}{4C_{s,q,\tilde q,d}},
$$
NSE has a mild solution $u$ on the interval $(0, T)$ so that 
\begin{equation}\label{3.35}
u \in \mathcal{K}_{q,\tilde q,T}^{s, \tilde q}.
\end{equation}
Lemma \ref{lem.3.4} and the relation \eqref{3.35} imply that  
$$
B(u,u) \in \mathcal{K}_{q,1,T}^{s,q} 
\subseteq \mathcal{K}_{q,r,T}^{s,q} \subseteq L^\infty\Big([0, T]; 
\dot{H}^s_{L^{q,r}}\Big).
$$ 
On the other hand, from Lemma \ref{lem.3.1}, we have 
$e^{t\Delta}u_0 \in L^\infty\Big([0, T]; 
\dot{H}^s_{L^{q,r}}\Big)$. \\
Therefore
$$
u = e^{t\Delta}u_0 - B(u,u) \in  L^\infty\Big([0, T]; \dot{H}^s_{L^{q,r}}\Big).
$$
From Lemma \ref{lem.3.5} and Lemma \ref{lem.3.3}, we deduce that $u \in \mathcal{K}_{q,1,T}^{s, \tilde q}$.\\
From the definition of $ \mathcal{K}_{q,r,T}^{s,\tilde q}$ and Lemma \ref{lem.3.5}, 
we deduce that the left-hand side of the inequality \eqref{3.34} converges to $0$ when $T$ 
tends to $0$. Therefore the inequality \eqref{3.34} holds for 
arbitrary $u_0 \in \dot{H}^s_{L^{q,r}}(\mathbb{R}^d)$ 
when $T(u_0)$ is small enough. \\
(b) From Lemma \ref{lem.2.1.4}, the two quantities 
$$
\big\|u_0\big\|_{\dot{B}^{\frac{d}{\tilde q} - 1, \infty}_{\tilde q}}
\ \ \text{and} \ \ 
\underset{0 < t < \infty}{\rm sup}t^{\frac{d}{2}
({\frac{1}{q} - \frac{1}{\tilde q}})}
\big\|e^{t\Delta}u_0\big\|_{\dot{H}^{\frac{d}{q} - 1}_{\tilde q}}
$$ 
are equivalent, then there exists a positive constant 
$\sigma_{q,\tilde q,d}$ such that if  \linebreak
$\big\|u_0\big\|_{\dot{B}^{\frac{d}{\tilde q} - 1, \infty}_{\tilde q}} 
\leq \sigma_{q,\tilde q,d}$ and  $T = \infty$ 
then the inequality \eqref{3.34} holds.\qed 
\begin{Not}
In the case when the initial data  belong to the critical Sobolev-Lorentz spaces 
$\dot{H}^{\frac{d}{q} - 1}_{L^{q,r}}(\mathbb{R}^d), (1<q \leq d, r \geq 1)$, from Theorem \ref{th.3.2} (b), 
we get the existence of global mild solutions in the spaces  
$L^\infty([0, \infty); \dot{H}^{\frac{d}{q} - 1}_{L^{q, r}}(\mathbb{R}^d))$  
when the norm of the initial value in the Besov spaces 
$\dot{B}^{\frac{d}{\tilde q} - 1, \infty}_{\tilde q}(\mathbb{R}^d)$ 
is small enough. Note that a function in $\dot{H}^{\frac{d}{q} - 1}_{L^{q, r}}(\mathbb{R}^d)$ 
can be arbitrarily large in the $\dot{H}^{\frac{d}{q} - 1}_{L^{q, r}}(\mathbb{R}^d)$ 
norm but small in the $\dot{B}^{\frac{d}{\tilde q} - 1, \infty}_{\tilde q}(\mathbb{R}^d)$ norm. 
This is deduced from the following imbedding maps (see Lemma \ref{lem.3.5})
$$
\dot{H}^{\frac{d}{q} - 1}_{L^{q, r}}(\mathbb{R}^d) \hookrightarrow 
 \dot{B}^{\frac{d}{\tilde q} - 1, \infty}_{\tilde q}(\mathbb{R}^d), 
\Big(\frac{1}{q}-\frac{1}{d} < \frac{1}{\tilde q} < \frac{1}{q}\Big).
$$ 
This result is stronger than that of Cannone. 
In particular, when $q = r=d, s=0$, we get back the Cannone theorem 
(Theorem 1.1 in \cite{M. Cannone 1997}).
\end{Not}
Next, we consider the super-critical indexes $s > \frac{d}{q} - 1$.
\begin{Th}\label{th.3.3.2}
Let
$$
s \geq 0, q > 1, r \geq 1,\ and\ \frac{s}{d} <  \frac{1}{q} 
< \frac{s + 1}{d}. 
$$
Then for any $\tilde q$ be such that
$$
\frac{1}{2}\Big({\frac{1}{q}+\frac{s}{d}}\Big) 
< \frac{1}{\tilde q} < {\rm min} \Big\{\frac{1}{2} 
+\frac{s}{2d}, \frac{1}{q}\Big\},
$$
there exists a positive constant $\delta_{s,q,\tilde q,d}$ 
such that for all $T > 0$ and for all \linebreak
 $u_0 \in \dot{H}^s_{L^{q,r}}(\mathbb{R}^d)\ 
with\ {\rm div}(u_0) = 0$ satisfying
\begin{equation*}
T^{\frac{1}{2}(1+s-\frac{d}{q})}\big\|u_0\big\|_{\dot{B}^{s- (\frac{d}{q} 
- \frac{d}{\tilde q}), \infty}_{\tilde q}} \leq \delta_{s,q,\tilde q,d},
\end{equation*}
NSE has a unique mild solution $u \in 
\mathcal{K}^{s, \tilde q}_{q,1,T} 
\cap L^\infty([0, T]; \dot{H}^s_{L^{q,r}}) $.
\end{Th}
\textbf{Proof}. Applying Lemma \ref{lem.2.1.4}, 
the two quantities $\big\|u_0\big\|_{\dot{B}^{s- (\frac{d}{q}
 - \frac{d}{\tilde q}), \infty}_{\tilde q}}$ 
and \linebreak $\underset{0 < t < \infty}
{\rm sup}t^{\frac{d}{2}(\frac{1}{q} - \frac{1}
{\tilde q})}\big\|e^{t\Delta}u_0\big\|_{\dot{H}^s_{\tilde q}}$ 
are equivalent. Thus
$$
\underset{0 < t < T}{\rm sup}t^{\frac{d}{2}(\frac{1}{q} 
- \frac{1}{\tilde q})}\big\|e^{t\Delta}u_0\big\|_{\dot{H}^s_{\tilde q}}
\lesssim \big\|u_0\big\|_{\dot{B}^{s- (\frac{d}{q} 
- \frac{d}{\tilde q}), \infty}_{\tilde q}},
$$
the theorem is proved by applying the above inequality 
and Theorem \ref{th.3.2}. \qed
\begin{Not}
In the case when the initial data  belong to the Sobolev-Lorentz spaces 
$\dot{H}^s_{L^{q,r}}(\mathbb{R}^d), (q > 1, r \geq 1, s \geq 0, 
\ \text{and}\ \frac{d}{q} - 1 < s  < \frac{d}{q})$, 
we obtain the existence of mild solutions in the spaces 
$L^\infty([0, T]; \dot{H}^s_{L^{q,r}}(\mathbb{R}^d))$ 
for any $T > 0$ when the norm of the initial value 
in the Besov spaces $\dot{B}^{s- (\frac{d}{q} 
- \frac{d}{\tilde q}), \infty}_{\tilde q}(\mathbb{R}^d)$ 
is small enough. Note that a function in $\dot{H}^s_{L^{q,r}}(\mathbb{R}^d)$ 
can be arbitrarily large in the $\dot{H}^s_{L^{q,r}}(\mathbb{R}^d)$ 
norm but small in $\dot{B}^{s- (\frac{d}{q} 
- \frac{d}{\tilde q}), \infty}_{\tilde q}(\mathbb{R}^d)$ norm. 
This is deduced from the following imbedding maps (see Lemma \ref{lem.3.5})
$$
\dot{H}^s_{L^{q,r}}(\mathbb{R}^d) \hookrightarrow   \dot{B}^{s- (\frac{d}{q} 
- \frac{d}{\tilde q}), \infty}_{\tilde q}(\mathbb{R}^d),\ \Big(\frac{s}{d} 
< \frac{1}{\tilde q} < \frac{1}{q}\Big).
$$
\end{Not}
Applying Theorem \ref{th.3.3.2} for $q > d, r=q$ and $s = 0$, 
we get the following proposition which is stronger than the 
result of Cannone and Meyer (\cite{M. Cannone 1995}, 
\cite{M. Cannone. Y. Meyer 1995}). 
In particular, we obtained a result that is stronger than that 
of Cannone and Meyer but under a much weaker condition 
on the initial data.
\begin{Md}\label{Md2}
Let $q > d$. Then for any $\tilde q$ be such that
$$
q < \tilde q < 2q,
$$
there exists a positive constant $\delta_{q,\tilde q,d}$ 
such that for all $T > 0$ and for all \linebreak
 $u_0 \in L^q(\mathbb{R}^d)\ with\ {\rm div}(u_0) = 0$ satisfying
\begin{equation}\label{th1}
T^{\frac{1}{2}(1-\frac{d}{q})}
\big\|u_0\big\|_{\dot{B}^{\frac{d}{\tilde q} 
- \frac{d}{q}, \infty}_{\tilde q}} \leq \delta_{q,\tilde q,d},
\end{equation}
NSE has a unique mild solution 
$u \in \mathcal{K}^{0, \tilde q}_{q,1,T} \cap L^\infty([0, T]; \L^q)$.
\end{Md}
\begin{Not}
If in \eqref{th1} we replace the $\dot{B}^{\frac{d}{\tilde q}
-\frac{d}{q}, \infty}_{\tilde q}$ norm by the $L^{q}$
norm then we get the assumption 
made in (\cite{M. Cannone 1995}, 
\cite{M. Cannone. Y. Meyer 1995}). We show that 
the condition \eqref{th1} is weaker 
than the condition in (\cite{M. Cannone 1995}, 
\cite{M. Cannone. Y. Meyer 1995}). 
In Remark 5 we have showed that
$$
L^q(\mathbb{R}^d) \hookrightarrow  \dot{B}^{\frac{d}{\tilde q}
-\frac{d}{q}, \infty}_{\tilde q}(\mathbb{R}^d), (\tilde q > q \geq d),
$$
but these two spaces are different. Indeed,  we have 
$\big|x\big|^{-\frac{d}{q}} \notin L^q(\mathbb{R}^d)$. 
On the other hand by using Lemma \ref{lem.2.1.4}, we can 
easily prove that
$\big|x\big|^{-\frac{d}{q}} \in \dot{B}^{\frac{d}{\tilde q}
-\frac{d}{q}, \infty}_{\tilde q}(\mathbb{R}^d)$ 
for all $\tilde q > q$.
\end{Not}
Applying Theorem \ref{th.3.3.2} for 
$q=r=2, \frac{d}{2} - 1 < s < \frac{d}{2}$, 
we get the following proposition which is stronger than the 
results of Chemin  
in \cite{J. M. Chemin 1992} and Cannone in \cite{M. Cannone 1995}. 
In particular, we obtained the result that is stronger than that 
of Chemin and Cannone but under a much weaker condition 
on the initial data.
\begin{Md}\label{Md3}
Let $\frac{d}{2} - 1 < s < \frac{d}{2}$. 
Then for any $\tilde q$ be such that
$$
\frac{1}{2}\Big({\frac{1}{2}+\frac{s}{d}}\Big) 
< \frac{1}{\tilde q} < \frac{1}{2}, 
$$
there exists a positive constant $\delta_{s,\tilde q,d}$ 
such that for all $T > 0$ and for all \linebreak
 $u_0 \in \dot{H}^s(\mathbb{R}^d)\ with\ {\rm div}(u_0) = 0$ 
satisfying
\begin{equation}\label{th3}
T^{\frac{1}{2}(1+s-\frac{d}{2})}
\big\|u_0\big\|_{\dot{B}^{s- (\frac{d}{2} 
- \frac{d}{\tilde q}), \infty}_{\tilde q}} \leq \delta_{s,\tilde q,d},
\end{equation}
NSE has a unique mild solution 
$u \in \mathcal{K}^{s, \tilde q}_{2,1,T} 
\cap L^\infty([0, T]; \dot{H}^s)$.
\end{Md}
\begin{Not}
If in \eqref{th3} we replace the $\dot{B}^{s- (\frac{d}{2} 
- \frac{d}{\tilde q}), \infty}_{\tilde q}$ norm by the $\dot{H}^s(\mathbb{R}^d) $
norm then we get the assumption 
made in (\cite{J. M. Chemin 1992}, \cite{M. Cannone 1995}). We show that 
the condition \eqref{th3} is weaker 
than the condition in (\cite{J. M. Chemin 1992}, \cite{M. Cannone 1995}). 
In Remark 5 we showed that
$$
\dot{H}^s(\mathbb{R}^d) \hookrightarrow   \dot{B}^{s- (\frac{d}{2} 
- \frac{d}{\tilde q}), \infty}_{\tilde q},\ \frac{1}{2}\Big({\frac{1}{2}
+\frac{s}{d}}\Big) < \frac{1}{\tilde q} < \frac{1}{2},
$$
but that these two spaces are different. Indeed, 
we have $\dot{\Lambda}^{-s}|.|^{-\frac{d}{2}} 
\notin \dot{H}^s(\mathbb{R}^d)$, on the other 
hand by using Lemma \ref{lem.2.1.4}, we easily prove that \linebreak
$\dot{\Lambda}^{-s}|.|^{-\frac{d}{2}} \in \dot{B}^{s- (\frac{d}{2} 
- \frac{d}{\tilde q}), \infty}_{\tilde q}(\mathbb{R}^d)$ for all $\tilde q  > 2$.
\end{Not}
\vskip 0.2cm 
{\bf Acknowledgments}. This research is funded by 
Vietnam National \linebreak Foundation for Science and Technology 
Development (NAFOSTED) under grant number  101.02-2014.50.

\end{document}